
\magnification=1200
\baselineskip=18pt
\def\I{I\!\!P}
\font\c=cmr8
\nopagenumbers
\headline={\ifnum\pageno=1 \else{\ifodd\pageno\rightheadline
\else\leftheadline\folio\fi}\fi}
\def\rightheadline{\c\hfil SOME
REMARKS ON THE OBSTRUCTEDNESS OF CONES\hfil\folio}
\def\leftheadline{\c\hfil CILIBERTO-LOPEZ-MIRANDA\hfil}
\voffset=1\baselineskip
\font\ms=msbm10
\font\msp=msbm5

\def\qed{\vrule height 1.1ex width 1.0ex depth -.1ex }
\def\mapright#1{\smash{\mathop{\longrightarrow}\limits^{#1}}}

\centerline{\bf SOME REMARKS ON THE OBSTRUCTEDNESS OF CONES }
\centerline{\bf OVER CURVES OF LOW GENUS}
\vskip 1.5cm
\centerline{CIRO CILIBERTO*\footnote\null{{\it \noindent 1991
Mathematics Subject Classification:} Primary 14C05. Secondary 14J10,
14B07.}, ANGELO FELICE LOPEZ\footnote*{Research partially supported by
the MURST national project ``Geometria Algebrica";}
AND RICK\footnote\null{the authors are members of GNSAGA of CNR.}
MIRANDA*\footnote*{* Research  partially supported by an NSA grant }}
\vskip 1.5cm \noindent {\bf 1. INTRODUCTION}
\vskip .5cm
Let $C$ be a smooth irreducible curve of genus $g$ and $L$ a very
ample line bundle on $C$ of degree $d$. Embedding $C$ in $\I^r =
\I H^0(L)^*$ we can consider the cone $X_{C,L} \subset \I^{r+1}$
over $C$ with vertex a point in $\I^{r+1} - \I^r$. When $g = 0$ it
was shown by Pinkham ([P]) that rational normal cones are obstructed,
that is they represent singular points of their Hilbert schemes, as
soon as $r \geq 4$. One of the purposes of the present article is to
employ the technique of Gaussian maps to study the obstructedness of
cones over curves of {\it positive} genus $g$ and degree $d >> g$.

For $k \geq 1$ set $R(\omega_{C} \otimes L^{k-1},L) = Ker
\{H^0(\omega_{C} \otimes L^{k-1}) \otimes H^0(L) \to H^0(\omega_{C}
\otimes L^k) \}$ and let   $$\Phi_{\omega_{C} \otimes L^{k-1},L} :
R(\omega_{C} \otimes L^{k-1},L) \to H^0(\omega_{C}^2 \otimes L^k)$$
\noindent be the Gaussian map defined locally by $\Phi_{\omega_{C}
\otimes L^{k-1},L}(s \otimes t) = s dt - t ds$. Moreover let
$\gamma_{C,L}^k = corank \ \Phi_{\omega_{C} \otimes L^{k-1},L}$, and
$\gamma_{C,L} = corank \ \Phi_{\omega_{C},L}$.

As is well-known now the integers $\gamma_{C,L}^k$ may be used to
compute the dimension of the tangent space to the Hilbert scheme at
points representing cones $X_{C,L}$ over $C$ (see for example [W1],
[St], [T], [CM], [CLM1], [CLM2]).

We will see in section 2 that there is an interesting difference for
the obstructedness of such cones between the nonhyperelliptic case
for $g \geq 3$ and the hyperelliptic case $g = 1, 2$. For $g \geq 3$
and $d >> g$ if $C$ is not hyperelliptic and $L$ is any line bundle on
$C$ of degree $d$,  we have $\gamma_{C,L} = 0$ ([W2], [BEL]), and the
cone $X_{C,L}$ is unobstructed. On the other hand if $C$ is
hyperelliptic and $d >> g$ then $\gamma_{C,L} > 0$ ([St]); in
particular if $g = 1, 2$ and $d \geq 3g + 7$, we will show in section
4, that the cone $X_{C,L}$ is obstructed and cannot be smoothed. By
taking general hypersurface sections of such cones and using the above
fact, we will construct infinitely many examples of {\it nonreduced}
components of the Hilbert scheme of curves in $\I^{r+1}, r \geq 9$
(see Theorem (4.11)). Note that so far all the known examples of
nonreduced components of the Hilbert scheme are for curves in $\I^3$
([M], [GP], [K], [E], [Fl]).

Another interesting feature of the coranks of Gaussian maps is that
they give the cohomology $h^0(N_{C}(-k)), k \geq 1$, of the normal
bundle of $C$ in $\I^r$.  This in turn, by a theorem of Zak [Z],
governs in many cases the existence of higher dimensional varieties
having $C$ as their curve section (see for example [CLM1], [CLM2] for
Fano varieties, [W1] for K3 surfaces). In section 3 we will further
exploit the power of the above technique, by applying it to varieties
of degree five having elliptic curve sections. In particular we will
extend to the singular case the well-known smooth classification
([Sc], [I], [Fu]).  We will then consider surfaces in $\I^r, r \geq 6$
with curve sections of genus at most three. In these cases a
straightforward application of the method of Gaussian maps and Zak's
theorem does not work. We will instead recover the known
classification by means of a projective technique (the ``tetragonal
lemma" (3.2)).

{\it Acknowledgments.} The authors wish to thank the Department of
Mathematics of the University of Trento and the CIRM for the warm
hospitality and the organization provided in the school-conference
``Higher dimensional complex varieties" held in Trento in June 1994.
We also thank M. Coppens for some helpful suggestions regarding lemma
(3.2).
\vskip .5cm
\noindent {\bf 2. GLOBAL SECTIONS OF THE NORMAL BUNDLE AND CORANK OF
GAUSSIAN  MAPS}
\vskip .5cm
As above we let $C$ be a smooth irreducible curve of genus $g \geq 0$,
$L$ be a line bundle on $C$ of degree $d \geq 1$ and $\gamma_{C,L} =
corank \ \Phi_{\omega_{C},L}$. We will first collect some facts about
$\gamma_{C,L}$.

{\bf Proposition (2.1).} {\sl We have the following values for
$\gamma_{C,L}$:

\noindent (2.2) If $g = 0$ then $\gamma_{C,L} = \cases{0 &{\sl for}
$d \leq 3$ \cr d-3 &{\sl for} $d \geq 4$ \cr}$;

\noindent (2.3) If $g = 1$ then $\gamma_{C,L} = d$;

\noindent (2.4) If $g = 2$ then $\gamma_{C,L} = 6$ for $d \geq 5$;

\noindent (2.5) If $g = 3$ and $C$ is not hyperelliptic we have
$$\gamma_{C,L} = \cases{0 &{\sl for} $d \geq 17$ \cr 0 &{\sl for $L$
general and} $d \geq 14$ \cr \geq 14 - d &{\sl for} $6 \leq d \leq 13$
\cr};$$

\noindent (2.6) If $g = 4$ and $C$ is not hyperelliptic we have
$$\gamma_{C,L} = \cases{0 &{\sl for} $d \geq 19$ \cr 0 &{\sl for $L$
general and} $d \geq 15$ \cr \geq 15 - d &{\sl for} $8 \leq d \leq 14$
\cr};$$

\noindent (2.7) If $g = 5$ and $C$ is not trigonal we have
$$\gamma_{C,L} = \cases{0 &{\sl for} $d \geq 17$ \cr 0 &{\sl for $L$
general and} $d \geq 12$ \cr \geq 36 - 3d &{\sl for} $10 \leq d \leq
11$ \cr};$$

\noindent (2.8) If $g \geq 6$ and both $C$ and $L$ are general we have
$$\gamma_{C,L} = 0 {\sl \ if \ } d \geq \cases{g+12 &{\sl for} $6 \leq
g \leq 8$ \cr g+9 &{\sl for} $g \geq 9$ \cr}.$$}
\noindent {\it Proof.} If $g = 0, 1$ we have $R(\omega_{C},L) = 0$,
hence $\gamma_{C,L} = h^0(\omega_{C}^2 \otimes L)$ and this gives
(2.2) and (2.3). If $g = 2$, by the base point free pencil trick, we
have $R(\omega_{C},L) = H^0(\omega_{C}^{-1} \otimes L)$ hence
$\gamma_{C,L} = h^0(\omega_{C}^2 \otimes L) - h^0(\omega_{C}^{-1}
\otimes L) + dim Ker \ \Phi_{\omega_{C},L} = d + 3 - ( d - 3 ) + dim
Ker \ \Phi_{\omega_{C},L} = 6 +  dim Ker \ \Phi_{\omega_{C},L}$. A
straightforward computation in the diagram
$$\matrix{H^0(\omega_{C}^{-1} \otimes L) \ \ \mapright{\Phi}  \ \
H^0(\omega_{C}^2 \otimes L) \cr \ \ \ \ \ \ \ \ \ \searrow \cong \ \ \
\ \ \ \ \nearrow \Phi_{\omega_{C},L} \cr R(\omega_{C},L) \cr}$$
\noindent gives that, if $H^0(\omega_{C}) = < s , t >$, we have
$\Phi(\sigma) = 2 \sigma \Phi_{\omega_{C},\omega_{C}}(s \otimes t)$,
and $\Phi_{\omega_{C},\omega_{C}}(s \otimes t) \not= 0$, hence both
$\Phi$ and $\Phi_{\omega_{C},L}$ are injective and (2.4) follows. To
see (2.5), (2.6) and (2.7) consider the canonical embedding $C \subset
\I^{g-1}$. As is well known (see [W2]), from the exact sequence
$$0 \to N^*_{C / \I^{g-1}} \otimes \omega_{C} \otimes L \to
\Omega_{{\I^{g-1}}_{|C}}^1 \otimes \omega_{C} \otimes L \to
\omega_{C}^2 \otimes L \to 0$$
\noindent we get

\noindent (2.9) \ \ \ $H^0(\Omega_{{\I^{g-1}}_{|C}}^1 \otimes
\omega_{C} \otimes L) \mapright{\Phi_{\omega_{C},L}} H^0(\omega_{C}^2
\otimes L) \to H^1(N^*_{C / \I^{g-1}} \otimes \omega_{C} \otimes L)
\to$

\noindent \ \ \ \ \ \ $\to H^1(\Omega_{{\I^{g-1}}_{|C}}^1 \otimes
\omega_{C} \otimes L)$

\noindent and in the cases at hand we will show that

\noindent (2.10) $H^1(\Omega_{{\I^{g-1}}_{|C}}^1 \otimes \omega_{C}
\otimes L) = 0$ and therefore

\noindent (2.11) $\gamma_{C,L} = h^1(N^*_{C / \I^{g-1}} \otimes
\omega_{C} \otimes L)$.

\noindent To see (2.10) consider the Euler sequence
$$0 \to \Omega_{{\I^{g-1}}_{|C}}^1 \otimes \omega_{C} \otimes L \to
H^0(\omega_{C}) \otimes L \to \omega_{C} \otimes L \to 0$$
\noindent and notice that, with the given hypotheses on $d$, we have
$H^1(L) = 0$ and the multiplication map $H^0(\omega_{C}) \otimes H^0(L)
\to H^0(\omega_{C} \otimes L)$ is surjective (by [C2]), hence (2.10).

When $g = 3$ by (2.11) we see that $\gamma_{C,L} = h^1(L \otimes
\omega_{C}^{-3}) = h^0(L^{-1} \otimes \omega_{C}^4)$ and this is zero
for $d \geq 17$ or is at least $14 - d$ for $6 \leq d \leq 13$. For
$14 \leq d \leq 16$ and $L$ general in $Pic^dC$ we have that
$\gamma_{C,L} = 0$ since $dim \ Pic^dC = 3$ while $dim \{L: h^0(L^{-1}
\otimes \omega_{C}^4) \not= 0 \} = 16 - d$. This gives (2.5). When
$g = 4$ we get by (2.11) that $\gamma_{C,L} = h^1(L \otimes
(\omega_{C}^{-1} \oplus \omega_{C}^{-2})) = h^0(L^{-1} \otimes
\omega_{C}^2) + h^0(L^{-1} \otimes \omega_{C}^3)$ and this is zero
for $d \geq 19$ or is at least $15 - d$ for $8 \leq d \leq 14$. For
$15 \leq d \leq 18$ and $L$ general in $Pic^dC$ we have that
$\gamma_{C,L} = 0$ since $dim Pic^dC = 4$ while $dim \{L: h^0(L^{-1}
\otimes  \omega_{C}^3) \not= 0 \} = 18 - d$. This proves (2.6).
Similarly when $g = 5$, again by (2.11), we have that  $\gamma_{C,L} =
h^1(L \otimes (\omega_{C}^{-1})^{\oplus 3}) = 3 h^0(L^{-1} \otimes
\omega_{C}^2)$ and this is zero for $d \geq 17$ or is at least $36 -
3d$ for $10 \leq d \leq 11$.  For $12 \leq d \leq 16$ and $L$ general
in $Pic^dC$ we have that $\gamma_{C,L} = 0$ since $dim Pic^dC = 5$
while $dim \{L: h^0(L^{-1} \otimes \omega_{C}^2) \not= 0 \} = 16 - d$.
Hence (2.7) follows. Finally (2.8) was proved in [L, Corollary (1.7)].
\ \qed

We turn now to the relation between the corank of Gaussian maps and
the obstructedness of cones. With notation as above, suppose from now
on that $L$ is very ample, $C \subset \I^r = \I H^0(L)^*$ and let
$X_{C,L}$ be a cone over $C$ in $\I^{r+1}$ with vertex a point $P \in
\I^{r+1} - \I^r$.

We recall for the reader's convenience the connection between
cohomology of the normal bundle and corank of Gaussian maps. Let
$N_{C}, N_{X_{C,L}}$ be the normal bundles of $C \subset \I^r$ and of
$X_{C,L} \subset \I^{r+1}$, respectively. Then

{\bf Proposition (2.12).} {\sl Suppose $g \geq 1$. Then:

\noindent (2.13) $h^0(N_{C}(-1)) = r + 1 + \gamma_{C,L}$;

\noindent (2.14) $h^0(N_{C}(-k)) = dim Coker \
\Phi_{\omega_{C} \otimes L^{k-1},L}$ for every $k \geq 2$;

\noindent (2.15) If $\gamma_{C,L} = 0$ then $h^0(N_{C}(-k)) = 0$ for
every $k \geq 2$;

\noindent (2.16) $h^0(N_{X_{C,L}}) \leq \sum\limits_{k \geq 0}
h^0(N_{C}(-k))$;

\noindent (2.17) If $C \subset \I^r = \I H^0(L)^*$ is projectively
normal we have equality in (2.16).}

\noindent {\it Proof.} (2.13) and (2.14) are in [CM, Proposition 1.2].
Of course (2.15) follows for $k \geq 3$ as soon as we prove it for $k =
2$. To this end from the diagram
$$\matrix{R(\omega_{C},L) \otimes H^0(L) &\mapright{\Phi_{\omega_{C},L}
\otimes id_{H^0(L)}} &H^0(\omega_{C}^2 \otimes L) \otimes H^0(L) \cr
\downarrow && \ \ \ \ \ \downarrow \mu \cr
R(\omega_{C} \otimes L,L) &\mapright{\Phi_{\omega_{C} \otimes L,L}}
&H^0(\omega_{C}^2 \otimes L^2) \cr}$$
\noindent we deduce that $\Phi_{\omega_{C} \otimes L,L}$ is surjective
since $\Phi_{\omega_{C},L}$ is by hypothesis and so is $\mu$ (see for
example [G, Theorem (4.e.1)]). Hence  $H^0(N_{C}(-2)) = 0$ by (2.14).

\noindent To see (2.16) set $X = X_{C,L}$ and notice that
${N_{X}}_{|C} \cong N_{C}$, and therefore the exact sequence
$$0 \to N_{X}(-h-1) \to N_{X}(-h) \to N_{C}(-h) \to 0$$
\noindent implies that $h^0(N_{X}(-h)) \leq h^0(N_{X}(-h-1))
+h^0(N_{C}(-h))$, and applying this successively we get (2.16) since
$h^0(N_{X}(-j)) = 0$ for $j >> 0$. Finally (2.17) follows by standard
facts (see for example [CM, (4.2)]). \ \qed

Now let $W$ be a component of the Hilbert scheme $H_{d,g,r}$ of curves
of degree $d$ and genus $g$ in $\I^r$ such that $W \owns [C]$, the
point representing $C$, and denote by ${\cal H}(W)$ the family of cones
in $\I^{r+1}$ over curves representing points of $W$. We have the
following general fact

{\bf Proposition (2.18).} {\sl Suppose that $g \geq 3$ and a general
point $[C_{\eta}]$ of $W$ is a smooth point such that
$\gamma_{C_{\eta},{\cal O}_{C_{\eta}}(1)} = 0$. Then

\noindent (a) ${\cal H}(W)$ is a generically smooth component of the
Hilbert scheme of surfaces of degree $d$ in $\I^{r+1}$ and $dim {\cal
H}(W) = h^0(N_{C_{\eta}}) + r + 1$. In particular $dim {\cal H}(W) =
(r + 1) (d + 1) + (r - 3) (1 - g)$ if $H^1(N_{C_{\eta}}) = 0$;

\noindent (b) If $\gamma_{C,L} = 0$ and $[C]$ is a smooth point of
$W$, then $X_{C,L}$ is unobstructed;

\noindent (c) If $C \subset \I^r = \I H^0(L)^*$ is projectively normal
and $X_{C,L}$ is unobstructed, then $\gamma_{C,L} = 0$.}

\noindent {\it Proof.} First suppose that $\gamma_{C,L} = 0$ and that
$[C]$ is a smooth point of $W$. Let $V$ be a component of the Hilbert
scheme of surfaces of degree $d$ in $\I^{r+1}$ such that $V \supseteq
{\cal H}(W)$. We have  $$h^0(N_{C}) + r + 1 = dim W + r + 1 = dim {\cal
H}(W) \leq dim V \leq  h^0(N_{X_{C,L}}) \leq h^0(N_{C}) + r + 1$$
\noindent where the last inequality follows by (2.16), (2.15) and
(2.13). Hence $V = {\cal H}(W)$ and $X_{C,L}$ is unobstructed. This
shows (b). Applying this to $C_{\eta}$ we get (a) since
$h^0(N_{C_{\eta}}) + r + 1 = (r + 1) (d + 1) + (r - 3) (1 - g)$ when
$H^1(N_{C_{\eta}}) = 0$. Now suppose $C \subset \I^r =  \I H^0(L)^*$
is projectively normal and $X_{C,L}$ is unobstructed. Since ${\cal
H}(W)$ is a component of the Hilbert scheme, it must be the only one
containing $[X_{C,L}]$. If $\gamma_{C,L} > 0$ we would have
$$\eqalign{& dim {\cal H}(W) = h^0(N_{C_{\eta}}) + r + 1 \leq
h^0(N_{C}) + r + 1 < \cr
& < h^0(N_{C}) + r + 1 + \gamma_{C,L} + \sum\limits_{k \geq 2}
h^0(N_{C}(-k)) = h^0(N_{X_{C,L}}) \cr}$$
\noindent by (2.17). But this shows that $X_{C,L}$ is obstructed.
Hence (c) is proved. \ \qed

{\bf (2.19)} {\it Remark.} Even for line bundles of high degree
it is possible that $\gamma_{C,L} > 0$ for a given pair $(C,L)$ while
$\gamma_{C_{\eta},L_{\eta}} = 0$ for the general pair
$(C_{\eta},L_{\eta})$. For example on a smooth nonhyperelliptic curve
$C$ of genus $g \geq 9$ take $P$ a general point of $C$ and $L =
\omega_{C}(2P)$. It is shown in [L] that $\gamma_{C,L} > 0$; on the
other hand by Proposition (2.1) we have $\gamma_{C_{\eta},L_{\eta}} =
0$ since  $deg  L_{\eta} = 2g \geq g + 9$. Other examples can be
deduced from [St], [T].

Pairing Propositions (2.1) and (2.18)  we get

{\bf Corollary (2.20).} {\sl Let $C$ be a smooth irreducible curve of
genus $g \geq 3$, $L$ a line bundle on $C$ of degree $d$ such that
either

\noindent (a) $d \geq 16 g - 2 g^2 - 13$ for $3 \leq g \leq 5$ or

\noindent (b) $d \geq 15 g - 2 g^2 - 13$ for $3 \leq g \leq 5$ and $L$
is general, $C$ is not hyperelliptic and also not trigonal if $g = 5$,
or

\noindent (c) $d \geq \cases{g+12 &{\sl for} $6 \leq g \leq 8$
\cr g+9 &{\sl for} $g \geq 9$ \cr}$ and both $C$ and $L$ are general.

Then $X_{C,L}$ is unobstructed and ${\cal H}(W)$ is a generically
smooth component of the Hilbert scheme of dimension $(r + 1) (d + 1) +
(r - 3) (1 - g)$ where $r = d - g$.}

\noindent {\it Proof.} By Proposition (2.1) under hypothesis (a), (b)
or (c), we have $\gamma_{C,L} = 0$. Moreover in all cases $H^1(L) = 0$,
hence we can apply (a) and (b) of Proposition (2.18). \ \qed
\vskip .5cm
\noindent {\bf 3. REMARKS ON THE CLASSIFICATION OF VARIETIES WITH
ELLIPTIC CURVE SECTIONS}
\vskip .5cm
Let $X^n \subset \I^N$ be a nondegenerate variety of
dimension $n \geq 2$ with $dim Sing X^n \leq n - 2$, degree $5$, whose
curve sections are smooth elliptic normal curves. If $H$ is the
hyperplane divisor of $X^n$ we have  $- K_{X^n} = (n - 1) H$ and
$N = n + 3$.

Well-known examples of such varieties are the linear sections of the
Grassmann variety {\ms G}$(1,4) \subset \I^9$ in its Pl\"ucker
embedding. Just as in [CLM2], we will show that the classification of
the varieties $X^n$ as above can be recovered by means of deformations
to cones and the theorem of Zak, whose statement we recall here.
A smooth nondegenerate variety $Y \subset \I^m$ is said to be
$k${\it-extendable} if there exists a variety $Z \subset \I^{m + k}$
that is not a cone and such that $Y = Z \cap \I^m$. Zak's theorem ([Z])
says that if $codim \ Y \geq 2$ and $h^0(N_{Y}(-1)) = m + 1$ then $Y$
is not 1-extendable; also if $k \geq 2, h^0(N_{Y}(-1)) \leq m + k$ and
$h^0(N_{Y}(-2)) = 0$, then $Y$ is not k-extendable.

In the Hilbert scheme of nondegenerate varieties $X^n \subset
\I^{n + 3}$ of dimension $n \geq 2$ with $dim Sing X^n \leq n - 2$,
degree 5 and such that $- K_{X^n} = (n - 1) H$, we let ${\cal X}_{n}$
be the open subset parametrizing varieties $X^n$ that are not cones
over an elliptic curve. Then

{\bf Proposition (3.1).} {\sl

\noindent (a) ${\cal X}_{n} \not= \emptyset$ if and only if $n \leq 6$;

\noindent (b) For $n \leq 6, {\cal X}_{n}$ is irreducible and the
family of linear  sections of {\ms G}$(1,4) \subset \I^9$ forms a dense
open subset of smooth points of ${\cal X}_{n}$.}

\noindent {\it Proof.} By (2.3) and (2.13) if $E \subset \I^4$ is a
smooth curve section of $X^n$, we have that $h^0(N_{E / \I^4}(-1)) =
10$. Also it is easily seen that $h^0(N_{E / \I^4}(-2)) = 0$, for
example using (2.14) and the surjectivity of $\Phi_{L,L}$ (see [BEL]).
Hence Zak's theorem (see [Z]) implies that $n \leq 6$.
On the other hand $dim \hbox{{\ms G}}(1,4) = 6$ and $-K_{\hbox{{\msp
G}}(1,4)} = 5H$, where $H$ is the Pl\"ucker divisor, hence
${\cal X}_{6} \not= \emptyset$ and so is ${\cal X}_{n}, n \leq 5$,
since a hyperplane section of an  $X^n$ is a $X^{n-1}$. This proves
(a). To see (b) we will use an argument similar to the one in [CLM2],
in the proof of theorems (3.2) and (3.11). Since the varieties $X^n, n
\geq 2$, are projectively Cohen-Macaulay (because $E$ is), they flatly
degenerate to the $n$-dimensional cone $\widehat E$ over their general
curve section $E$. Of course the locus of such cones is irreducible,
hence (b) follows if we show that these cones are smooth points of the
closure of ${\cal X}_{n}$. To this end we will prove, as in [CLM2],
that $h^0(N_{\widehat E})$ is bounded above by the dimension of the
family of the known examples, that is linear sections of {\ms G}$(1,4)
\subset \I^9$. Since $h^0(N_{E}(-k)) = 0$ for every $k \geq 2$, we
have  $$h^0(N_{\widehat E}) = h^0(N_{E}) + (n-1) h^0(N_{E}(-1)) = 25 +
10(n-1).$$

\noindent On the other hand the family of
$X^n \subset \I^{n + 3}$ obtained by linear sections of {\ms G}$(1,4)
\subset \I^9$ has dimension

\noindent $dim \hbox{{\ms G}}(n+3,9) + dim Aut \I^{n + 3} -
dim \hbox{Aut{\ms G}}(1,4) = (n+4)(6-n) + (n+4)^2 - 25 = 15 + 10n.$
\ \qed

Note that the case $n = 2$ of Proposition (3.1) is the one of Del Pezzo
surfaces of  degree $5$ in $\I^5$. Suppose now $S \subset \I^r, r \geq
6$ is a (not necessarily smooth) surface whose general hyperplane
section is a smooth elliptic normal curve $E$. We have  $deg E = r$,
hence $S$ is a surface of degree $r$ in $\I^r$ nonsingular in
codimension one. In this case a straightforward application of the
methods of [CLM2] does not allow to recover the family of Del Pezzo
surfaces. In fact by (2.3) and (2.13) we have $h^0(N_{E}(-1)) = 2r$,
hence Zak's theorem does not rule out the existence of such surfaces
$S \subset \I^r$ for $r \geq 10$. Moreover even for $6 \leq r \leq 9$
the upper bound on $h^0(N_{\widehat E})$ is larger than the dimension
of the family of Del Pezzo surfaces. Also note that this implies that
the cones over elliptic normal curves are obstructed for $6 \leq r
\leq 9$; in section 4 we will see that they are also obstructed for $r
\geq 10$.

We will then use another strategy to recover the classification.
Recall that $E$ is projectively Cohen-Macaulay, its homogeneous ideal
is generated by quadrics, the relations among them being minimally
generated by linear ones. The same then holds for $S$. Let $C = S \cap
Q$ be a smooth quadric section of $S$. We have that $C$ is a canonical
curve whose  homogeneous ideal is generated by quadrics and therefore
by Petri's theorem $C$ is not trigonal nor isomorphic to a plane
quintic. Hence $Cliff(C) \geq 2$, where $Cliff(C)$ is the Clifford
index of $C$. Moreover the Koszul relations among $Q$ and the quadrics
generating the ideal of $S$ give rise to nonlinear syzygies among the
generators of the ideal of $C$, relations that do not depend on the
linear ones (or in other words $h^0(N_{C}(-2)) \not= 0$). By a result
of Schreyer [S1] and Voisin [V], we must have $Cliff(C) = 2$, that is
$C$ is either tetragonal or isomorphic to a plane sextic.

This fact will allow us to classify $S$. The key point is given by the
following

{\bf Lemma (3.2)} (tetragonal lemma). {\sl Let $S \subset \I^r, r \geq
5$ be a surface, nonsingular in codimension one, not a cone, whose
general hyperplane  section has genus $g$, degree $d \geq 2g+3$ and is
linearly normal. If there  is a linear system of curves on $S$, of
dimension at least $1$, whose general element $D$ is smooth,
irreducible, special and tetragonal, then $S$ is either a rational
normal scroll (i.e. $g = 0$) or it lies on a rational normal threefold
scroll of planes, each meeting $S$ in a conic.}

\noindent {\it Proof.} First notice that $S$ is cut out by quadrics
because its general hyperplane section is projectively normal and is
cut out by  quadrics, since $d \geq 2g + 3$ [G]. Take $p_{1} + \ldots +
p_{4}$ a general divisor of the $g_{4}^1$ on $D$. The span of this
divisor in $\I^r$ is at most a 2-plane $\pi$. If it is a line, then $S$
is a rational normal scroll. Suppose now $\pi$ is a 2-plane. Consider
the hyperplane sections $\{ H_{t} \}$ of $S$ through $\pi$. Suppose
that a general $H_{t}$ is irreducible. By monodromy the behaviour of
$H_{t}$ at the points $p_{1}, \ldots , p_{4}$ is the same, and we claim
that $p_{1}, \ldots , p_{4}$ are smooth points of $H_{t}$. Otherwise
$H_{t}$ would be tangent at $p_{1}, \ldots , p_{4}$ and therefore $\pi$
would be tangent to $S$ at these four points. Now let $D$ and $p_1$
vary. Since $p_{1}$ describes the whole surface, we conclude that the
general tangent plane to the surface is tangent at four points, a
contradiction, since the Gauss map is birational. By the hypothesis $d
\geq 2g+3$, if $H_{t}$ is irreducible, we have that $H^1({\cal
O}_{H_{t}}(1)(- p_{1} - \ldots - p_{4})) = 0$ and hence $dim  |{\cal
O}_{H_{t}}(1)(- p_{1} - \ldots - p_{4})| = d - 4 - g$ while the
$H_{t'}$ cut out, away from $p_{1}, \ldots , p_{4}$, a linear series on
$H_{t}$ of dimension at least $r - 3 = d - g - 3$. This contradiction
shows that a general $H_{t}$ must be reducible. Therefore, by Bertini's
theorem, there is a fixed component $F$ of the linear system $\{ H_{t}
\}$, and $F \subset \pi$ unless  $\{ H_{t} \}$ is composed with a
pencil. In the latter case, since the dimension of $\{ H_{t} \}$ is at
least $r-3 \geq 2$, we have that $H_t$ has to contain the tangent plane
to $S$ at each point $p_{i}, i = 1, \ldots , 4$, for all $t$'s and
this, as we saw, is impossible. Since $S$ is cut out by quadrics, then
the degree of $F$ is at most $2$ and $F$ has to contain  $p_{1},\ldots
,p_{4}$. If $F$ is a line we conclude that $S$ is a rational normal
scroll. If $F$ is a conic, and is reducible, then $S$ is a scroll.
Arguing as we did above, we see that the lines of the scroll are
pairwise coplanar, hence $S$ is a cone, which we excluded. If $F$ is an
irreducible conic, then $S$ is contained in the rational normal
threefold scroll of planes $\bigcup\limits_{p_{1} + \ldots + p_{4} \in
g_{4}^1} <p_{1} + \ldots + p_{4}>$. \ \qed

We will use the tetragonal lemma to recover the classification of Del
Pezzo surfaces.

{\bf Proposition (3.3).} {\sl Let $S \subset \I^r, r \geq 6$ be a
surface whose general hyperplane section is a smooth elliptic normal
curve $E$. Then $S$ is either an elliptic normal cone or $r \leq 8$ and
$S$ is a divisor of class $2 H + (4 - r) R$ on a rational normal
threefold scroll of planes $X \subset \I^r$, with hyperplane section
$H$ and ruling $R$, or $S$ is the Veronese embedding $v_{3}(\I^2)
\subset \I^9$ given by the cubics.}

\noindent {\it Proof.} Notice that if $S$ is a cone then it is an
elliptic normal cone. From now on we suppose that $S$ is not a cone. We
claim that either $S$ lies on a rational normal threefold scroll of
planes, each meeting $S$ in a conic, or $r = 9$ and $S$ is smooth. In
fact if $S$ has a singular point $P$, since the ideal of $S$ is
generated by quadrics, the projection $S' \subset \I^{r-1}$ of $S$ from
$P$ is a surface of degree at most $r - 2$, hence $deg  S' = r - 2$ and
$S'$ contains a one dimensional family of lines $\{ L_{t} \}$. Set
$\pi_{t} = <P, L_{t}>$ and $F_{t} = \pi_{t} \cap S$. Then $deg  F_{t}
\leq 2$ since the ideal of $S$ is generated by quadrics and $S$ lies on
a rational normal threefold scroll of planes, each meeting $S$ in a
conic. If $S$ is nonsingular, by the discussion at the beginning of
this section, we know that a smooth quadric section $C$ of $S$ is
either tetragonal or isomorphic to a plane sextic. In the former case
we can apply the tetragonal lemma since $H^1({\cal O}_{C}(1)) =
H^1(\omega_{C}) \not= 0$, concluding again that $S$ lies on a rational
normal threefold scroll of planes, each meeting $S$ in a conic.  When
$S$ lies on a scroll as above, pulling back to the normalization, we
have $S \sim a H + b R$ and $2 = deg R \cap S = R \cdot S \cdot H = R
\cdot (a H^2 + b H \cdot R) = a$ since $R^2 = 0, H^2 \cdot R = 1$.
Moreover $r = deg S = S \cdot H^2 = 2(r-2) + b$, hence $b = 4 - r$. By
[S2, (6.3)] we have $2 deg X + 3 b \geq 0$, that is $r \leq 8$. It
remains the case where smooth quadric sections $C$ of $S$ are
isomorphic  to a plane sextic. Then $r = 9$ and by what we just showed,
$S$ must be smooth. Moreover $K_{S} = - H, P_{2}(S) = h^0({\cal
O}_{S}(-2)) = 0$ and $q(S) = h^1({\cal O}_{S}(-1)) = 0$, therefore $S$
is rational. Also $S$ is minimal, since by Noether's formula we get
$b_{2}(S) = h^{1,1}(S) = 1$, hence $S = v_{3}(\I^2)$. \ \qed

{\bf (3.4)} {\it Remarks.}

\noindent {\it (i)} With the tetragonal lemma it is easily recovered
also the classification of surfaces with sectional genus $2$.  Let $S
\subset \I^r, r \geq 6$ be a surface whose general hyperplane section
is a smooth linearly normal curve $C_{1}$ of genus  $2$. Then $S$ is
either a cone over $C_{1}$, or $r \leq 11$ and $S$ is a divisor of
class  $2 H + (5 - r) R$ on a rational normal threefold scroll of
planes $X \subset \I^r$, with hyperplane section $H$ and ruling $R$.

\noindent {\it Proof.} Let $C = S \cap Q$ be a smooth quadric section
of $S$. We have  $deg C = 2r + 2, g(C) = r + 4$ and $C$ is linearly
normal and special since  $h^1({\cal O}_{C}(1)) = 2$. Moreover $deg
\omega_{C}(-1) = 4$ and therefore $S$ has a linear system of dimension
at least $2$, whose general element is smooth, irreducible, special and
tetragonal (of course $C$ is not trigonal, since its ideal is generated
by quadrics). Since $deg  S =  r + 1 \geq 7$, by Lemma (3.2) we have
that, if $S$ is not a cone, then $S \sim a H + b R$ on a rational
normal threefold scroll of planes $X \subset \I^r$ of degree $r - 2$,
and each plane meets $S$ in a conic. Therefore, as in the proof of
Proposition (3.3), we get $a = 2$ and $b = 5 - r$. By [S2, (6.3)] we
have $2deg X + 3b \geq 0$, that is $r \leq 11$. \ \qed

\noindent {\it (ii)} Now let $S \subset \I^r, r \geq 6$ be a surface
whose general hyperplane section is a smooth linearly normal curve
$C_{1}$ of genus  $3$. Then $S$ is either a cone over $C_{1}$, or $r
\leq 14$ and $S$ is the Veronese embedding $v_{4}(\I^2) \subset
\I^{14}$ given by the quartics or a projection of it in $\I^{14-r}$ or
$S$ is a quadric section of the cone in $\I^6$ over the Veronese
surface  $v_{2}(\I^2) \subset \I^5$. The proof of this is similar to
the one above and will be omitted.

\noindent {\it (iii)} Suppose now $S \subset \I^r, r \geq 4$, is a
surface of degree $d$ whose general hyperplane section is a smooth
linearly normal curve $C$ of genus  $g \geq 3$. As already observed by
many authors ([BEL], [BC]), by using Zak's theorem ([Z]), one can
obtain an upper bound for $d$. If $C$ is hyperelliptic or trigonal this
is well known (see [Se1], [Se2], [Fa]). If $Cliff(C) = 2$ then $d \leq
4g - 4$, while if $Cliff(C) \geq 3$ then $d \leq 4g - 3 Cliff(C)$. The
above bounds follow by [BEL], since, if $d$ is larger, then
$\gamma_{C,{\cal O}_{C}(1)} = 0$.

\vskip .5cm
\noindent {\bf 4. NONREDUCED COMPONENTS OF THE HILBERT SCHEME OF CUR-}

\noindent {\bf VES}
\vskip .5cm
Let $C$ be a smooth irreducible curve of genus $g \geq 1$ and
$L$ a line bundle on $C$ of degree $d$ such that $d \geq 10$ for
$g = 1, d \geq 4g + 5$ for $g \geq 2$. As opposite to the non
hyperelliptic case $g \geq 3$, we saw in Proposition (2.1) that in case
$g = 1, 2$, we have $\gamma_{C,L} > 0$. Even though Proposition (2.18)
does not apply, we will see here that in fact the cone $X = X_{C,L}$
over $C$ with vertex a point $P \in \I^{r+1} - \I^r$ is obstructed if
$g = 1, 2$.

Note that $L$ is very ample, nonspecial, and, by well-known results,
$C \subset \I^r = \I H^0(L)^*$ is projectively normal, its ideal is
generated by quadrics and the relations among them are linear ([G]).
Now consider a general hypersurface $F_{n}$ in $\I^{r+1}$ of degree $n
\geq 4$ and let $\Gamma^n = X \cap F_{n}$. The nonreduced components of
the Hilbert scheme will be obtained by the family of such $\Gamma^n$.

Let us first collect some information on $X$ and $\Gamma^n$.

{\bf Proposition (4.1).} {\sl

\noindent (4.2) $\Gamma^n$ is a projectively normal smooth curve of
degree $nd$ and genus $p_{n} = ng + d {n (n-1) \over 2} - n + 1$;

\noindent (4.3) Let ${\cal F} = \{ \Gamma_{t} \}_{t \in D}$ be a flat
family of projective curves parametrized by a smooth irreducible
variety $D$ such that:

(i) ${\cal F}$ is a projective family of curves in $\I^{r+1}$;

(ii) there exists a closed point $t_{0}$ in $D$ such that
$\Gamma_{t_{0}}$ is a smooth curve of degree $nd$ and genus $p_{n} = ng
+ d {n (n-1) \over 2} - n + 1$, complete intersection of a hypersurface
of degree $n \geq 4$ with a cone over a projectively normal, nonspecial
curve $C$ of degree $d$ and genus $g$.

Then there is a neighborhood $U$ of $t_{0}$ in $D$ in the Zariski
topology such that, for all closed points $t \in U$, $\Gamma_{t}$ is
again the complete intersection of a hypersurface of degree $n$ with a
cone over a projectively normal, nonspecial curve of degree $d$ and
genus $g$.}

\noindent {\it Proof.} 	Let $\Delta$ be a general hyperplane section
of $\Gamma^n$; we may think of $\Delta$ as a divisor on $C$ belonging
to the linear system $|{\cal O}_{C}(n)|$. If we denote by
$h_{\Delta}(t)$ the Hilbert function of $\Delta$, we find
$$h_{\Delta}(t) = \cases{h^0(C,{\cal O}_{C}(t)) = td - g + 1, &
\hbox{for} t = 1,..., n-1 \cr h^0(C,{\cal O}_{C}(n)) - 1 = nd - g &
\hbox{for} t = n \cr nd & \hbox{for}  $t \geq n+1$ \cr}.$$ \noindent
Therefore we have $$\sum\limits_{t = 1}^n ( nd - h_{\Delta}(t)) =
\sum\limits_{t = 1}^n  [g + (n - t) d - 1] + 1 = ng + d {n (n-1) \over
2} - n + 1 = p_{n}$$ \noindent which, in view of [C1], remark (1.8),
(ii), implies that $\Gamma^n$  is projectively normal. To see (4.3) we
can assume $D$ to be affine, thus $D = Spec (A)$. We have then a family
$Y \subset \I_{A}^{r+1} = Proj A[x_{0},...,x_{r+1}]$ over $Spec (A)$.
If ${\cal I}_{Y}$ is the ideal sheaf  of $Y$ in $\I^{r+1}$, then ${\cal
I}(Y) = \bigoplus\limits_{n \in \hbox{\ms N}} H^0({\cal I}_{Y}(n))$ is
the homogeneous ideal of $Y$ in $S = A[x_{0},...,x_{r+1}]$. Consider
the homogeneous ideal ${\cal I} \subset {\cal I}(Y)$ generated by
${\cal I}(Y)_{2} = H^0({\cal I}_{Y}(2))$, and look at the scheme $W =
Proj (S / {\cal I})$. There is a commutative diagram $$\matrix{& Y
\hookrightarrow  W \subset & \I_{A}^{r+1} = Proj  A[x_{0},...,x_{r+1}]
\cr & \searrow \ \ \downarrow \cr & \ \ Spec \ (A) \cr}$$
\noindent and restricting to the local ring of $D$ at $t_{0}$ we find
that the scheme $W$ is flat over the generic point $\xi$ of $D$ (the
proof of this is similar to the one of proposition (1.6) of [C1] and
will be omitted). Therefore $W$ is flat over a neighborhood $U$ of
$t_{0}$ in $D$, thus giving a flat family of surfaces $X_{t}$, with
central fiber the cone $X$.  Notice that since $C$ is smooth and
projectively normal, then $X$ is normal  at the vertex $P$. For a
general point $\xi$ of $D, X_{\xi}$ has to be singular (see [P, theorem
(7.5)]), and has normal singularities (see [EGA, theorem 12.2.4]). On
the other hand, by theorem 2, chapter 1 of Horowitz's thesis [Ho],
$X_{\xi}$ (which is irreducible and nondegenerate, as well as $X$) is
a  scroll. Therefore $X_{\xi}$ has to be a cone, because these are the
only scrolls which may have normal singularities as we will see in the
Claim (4.4) below. In any case, since the family $\{ \Gamma_{t} \}$ is
very flat, lifting the equation of the hypersurface $F_{n}$, we see
that also $\Gamma_{t}$ is the complete intersection of a hypersurface
of degree $n$ with a cone over a projectively normal, nonspecial curve
of degree $d$ and genus $g$.

{\it Claim} (4.4). {\sl Let $T$ be an irreducible, nondegenerate normal
surface in $\I^r, r \geq 3$ and assume that $T$ is a scroll (that is
$T$ is ruled by  lines). If $T$ is singular then $T$ is a cone over a
smooth curve.}

\noindent {\it Proof.} Let $Q$ be a singular point of $T$ and let $H$
be the hyperplane section of $T$ with a general hyperplane passing
through $Q$. We shall  separately discuss two cases.

{\it Case 1:} $H$ is reducible.

\noindent Let $\pi : {\widetilde T} \to T$ be the minimal
desingularization of $T$ and let $E$ be the divisor contracted by $\pi$
in $Q$. Let moreover $\Sigma$ be the linear system of divisors on
$\widetilde T$ which are the pull-backs on  $\widetilde T$ of the
hyperplanes of $\I^r$ passing through $Q$. Clearly the  base locus of
$\Sigma$ coincides with $E$ and  our hypothesis on $H$ yields that the
movable part of $\Sigma$ is composed with a pencil $\cal P$.  Now let
$R$ be a general point of $T$ and $\widetilde R$ the corresponding
point of $\widetilde T$. All curves of $\Sigma$ passing through
$\widetilde R$ must contain the curve of $\cal P$ passing through
$\widetilde R$. This curve, in turn, has to be the pull-back on
$\widetilde T$ of some curve on $T$ contained in all the hyperplanes
passing through $Q$ and $R$. Whence clearly the general curve in $\cal
P$ is the pull-back on $\widetilde T$ of a line contained in $T$ and
passing through $Q$. This implies that $T$ is a cone with $Q$ as
vertex. The directrix of the cone  has to be smooth by the normality of
$T$.

{\it Case 2:} $H$ is irreducible.

\noindent It is easy to construct in this case a flat family
$\{ H_{t} \}_{t \in B}$ such that:

\noindent (i) $B$ is a smooth curve;

\noindent (ii) for every $t \in B, H_{t}$ is a hyperplane section of
$T$;

\noindent (iii) there exists a closed point $t_{1} \in B$ such that for
all $t \in B - \{ t_{1} \}, H_{t}$ is a smooth irreducible hyperplane
section of $T$;

\noindent (iv) in correspondence with the point $t_{1} \in B$ one has
that $H_{t_{1}}$ is a general hyperplane section of $T$ through $Q$;

\noindent (v) $H_{t_{1}}$ is reduced (because of the $S_{2}$ property
of $T$ at $Q$) irreducible and singular only at $Q$.

In view of our hypothesis on $H$, the lines of the scroll do not all
pass through $Q$ (otherwise $T$ would be a cone and $H$ reducible). But
this clearly implies that $H_{t}$, for $t \not= t_{1}$, is birational
to $H_{t_{1}}$, a  contradiction since the geometric genus of
$H_{t_{1}}$ has to be strictly less  than the one of $H_{t}$ (see
[Hi]).

\noindent This proves Claim (4.4) and hence also concludes the proof of
Proposition (4.1). \ \qed

Now let $N_{X}$ and $N_{\Gamma^n}$ be the normal bundles of $X$ and
$\Gamma^n$, respectively in $\I^{r+1}; N_{\Gamma^n / X}$ the normal
bundle of $\Gamma^n$ in $X$. We have

{\bf Proposition (4.5).} {\sl Let $C \subset \I^r = \I H^0(L)^*$ be a
smooth irreducible curve of genus $g \geq 1, L$ a line bundle on $C$ of
degree $d \geq \cases{10 & \hbox{for} $g = 1$ \cr 4g + 5 & \hbox{for}
$g \geq 2$ \cr}$ and let $X$ be the cone over $C$ with vertex a point
$P \in \I^{r+1} - \I^r$. Then

\noindent (4.6) $H^0(X,N_{X}(-k)) = 0$ for all $k \geq 2$;

\noindent (4.7) $H^1(\Gamma^n, N_{\Gamma^n / X}) = 0$;

\noindent (4.8) $h^0(\Gamma^n, N_{\Gamma^n}) = h^0(X,N_{X}) +
h^0(\Gamma^n, N_{\Gamma^n / X})$.}

\noindent {\it Proof.} To see (4.6) first notice that by standard facts
(for example applying (2.14) and [BEL]) it follows that

\noindent (4.9) $H^0(C,N_{C}(-k)) = 0$ for all $k \geq 2$.

Of course the assertion (4.6) is true for $k >> 0$. So we may argue by
descending induction. From the exact sequence
$$0 \to N_{X}(-k-1) \to N_{X}(-k) \to N_{C}(-k) \to 0$$
\noindent we see that (4.6) follows by (4.9). Since $N_{\Gamma^n / X}
\cong {\cal O}_{\Gamma^n}(n)$, we have
$$deg N_{\Gamma^n / X} = n^2 d > n^2 d - n (d - 2g + 2)	=
2gn + dn(n-1) - 2n = 2p_{n} - 2$$
\noindent hence (4.7). To show (4.8) observe that by (4.6) we have
$H^0(X,N_{X}(-n)) = 0$, hence there is an injection
$H^0(X,N_{X}) \mapright{h} H^0(\Gamma^n,N_{{X}_{|\Gamma^n}})$. Now
consider the exact sequence
$$0 \to H^0(\Gamma^n, N_{\Gamma^n / X}) \to H^0(\Gamma^n, N_{\Gamma^n})
\mapright{f} H^0(\Gamma^n,N_{{X}_{|\Gamma^n}}) \to 0 \leqno (4.10)$$
\noindent the map $f$ being surjective by (4.7). By the proof of (4.3)
we can deduce the existence of a map $\phi$ which makes the following
diagram commutative
$$\matrix{H^0(\Gamma^n, N_{\Gamma^n})
\ \mapright{f} \ H^0(\Gamma^n,N_{{X}_{|\Gamma^n}}) \cr
\ \ \ \ \ \ \searrow \phi \ \ \ \ \ \ \nearrow h \cr
\ \ \ \ \ H^0(X,N_{X}) \cr}$$
\noindent as one can see by a local computation. This implies that $h$
is also surjective. The assertion then follows by (4.10). \ \qed

We shall finally prove the announced result on nonreduced components
of the Hilbert scheme.

{\bf Theorem (4.11).} {\sl Let $\Gamma^n$ be a smooth complete
intersection with a hypersurface of degree $n \geq 4$ of a cone $X =
X_{C,L}$ over a smooth  irreducible curve $C$ of genus $g = 1, 2$,
embedded in $\I^r = \I H^0(L)^*$ by a line bundle $L$ of degree $d \geq
3g + 7$. Then $\Gamma^n$ belongs to a unique irreducible component
$\cal V$ of the Hilbert scheme $H_{nd,p_{n},r+1}, r = d - g$, such that
$\cal V$ is nonreduced and  $$dim {\cal V} = (r + 1)(d + 1) + (r - 3 +
n)(1 - g) + d{n(n+1) \over 2}.$$} \noindent {\it Proof.} Let $\cal V$
be a component of the Hilbert scheme containing the  point $[\Gamma^n]$
and $W$ the component of the Hilbert scheme containing $[C]$. The proof
of (4.3) yields the existence of a morphism  $\psi : {\cal V}_{red} \to
{\cal H}(W)_{red}$, which is clearly surjective, its (set-theoretical)
general fiber being nothing else than the projective space
corresponding to the linear system of curves cut out on a cone by the
hypersurfaces of degree $n$ in $\I^{r+1}$. Since for any such smooth
curve $\Gamma'$ on a cone $X'$ we have, as we saw in (4.7), that
$H^1(\Gamma', N_{\Gamma' / X'}) = 0$, then  $h^0(\Gamma', N_{\Gamma' /
X'}) = n^2 d - (ng + d {n(n-1) \over 2} - n + 1) + 1 = d {n(n+1) \over
2} - (g - 1) n$. Therefore we find
$$\eqalign{dim {\cal V} & = dim {\cal H}(W) + d {n(n+1) \over 2} - (g
- 1) n = dim W + r + 1 + d {n(n+1) \over 2} - (g - 1) n = \cr
& = (r + 1)(d + 1) + (r - 3 + n)(1 - g ) + d {n(n + 1) \over 2}. \cr}$$
\noindent Thus the part of the statement concerning the uniqueness of
${\cal V}$ and its dimension is proved. Now we notice that, by
Proposition (2.1) and by (2.13) and (4.8), we have
$$\eqalign{& dim {\cal V} = h^0(\Gamma^n, N_{\Gamma^n / X}) + dim {\cal
 H}(W) = h^0(\Gamma^n, N_{\Gamma^n / X}) + h^0(C, N_{C}) + r + 1 < \cr
& < h^0(\Gamma^n, N_{\Gamma^n / X}) + h^0(C, N_{C}) + r + 1 +
\gamma_{C,L} \leq h^0(\Gamma^n, N_{\Gamma^n / X}) + h^0(X,N_{X}) = \cr
& = h^0(\Gamma^n, N_{\Gamma^n}) = dim T_{[\Gamma^n]}{\cal V} \cr}$$
where $T_{[\Gamma^n]}{\cal V}$ is the Zariski tangent space to
${\cal V}$ at $[\Gamma^n]$. \ \qed

{\bf (4.12)} {\it Remark.} We limited ourselves to constructing
nonreduced components of the Hilbert scheme of curves that are complete
intersection of cones with hypersurfaces. Similar results can be proved
also for all curves of  sufficiently high degree on such cones and for
curves which are algebraically equivalent to a hypersurface
intersection plus a line.

\vskip .5cm
\noindent {\bf REFERENCES}
\baselineskip 12pt
\vskip .5cm
\item{[BC]} Ballico,E., Ciliberto,C.:\ On gaussian maps for projective
varieties.\ In:\ {\it Geometry of complex projective varieties, Cetraro
(Italy), June 1990.\ Seminars and Conferences \bf 9}.\ Mediterranean
Press:\ 1993, 35-54.
\vskip .2cm
\item{[BEL]} Bertram,A., Ein,L., Lazarsfeld,R.:\ Surjectivity of
Gaussian maps for line bundles of large degree on curves.\ In:\ {\it
Algebraic Geometry, Proceedings Chicago 1989.\ Lecture Notes in Math.\
\bf 1479}.\ Springer, Berlin-New York:\ 1991, 15-25.
\vskip .2cm
\item{[C1]} Ciliberto,C.:\ On the Hilbert scheme of curves of maximal
genus in a projective space.\ {\it Math.\ Z.\ \bf 194}, (1987),
351-363.
\vskip .2cm
\item{[C2]} Ciliberto,C.:\ Sul grado dei generatori dell'anello
canonico di una superficie di tipo generale.\ {\it Rend.\ Sem.\ Mat.\
Univ.\ Politecn.\ Torino \bf 41}, n.\ 3 (1983), 83-111.
\vskip .2cm
\item{[CLM1]} Ciliberto,C., Lopez,A.F., Miranda,R.:\ Projective
degenerations of K3 surfaces, Gaussian maps and Fano threefolds.\
{\it Invent.\ Math.\ \bf 114}, (1993) 641-667.
\vskip .2cm
\item{[CLM2]} Ciliberto,C., Lopez,A.F., Miranda,R.:\ Classification of
varieties with canonical curve section via Gaussian maps on canonical
curves.\ {\it Preprint} (1994).
\vskip .2cm
\item{[CM]} Ciliberto,C., Miranda,R.:\ On the Gaussian map for
canonical curves of low genus.\ {\it Duke Math.\ J.\ \bf 61}, (1990)
417-443.
\vskip .2cm
\item{[E]} Ellia,P.:\ D'autres composantes non r\`eduites de
$Hilb \I^3$.\ {\it Math.\ Ann.\ \bf 277}, (1987) 433-446.
\vskip .2cm
\item{[EGA]} Grothendieck,A., Dieudonn\`e,J.:\ EGA IV, Etude locale des
sch\'emas et des morphismes de sch\'emas.\ {\it Publ.\ Math.\ IHES \bf
28}, (1966).
\vskip .2cm
\item{[Fa]} Fania,M.L.:\ Trigonal hyperplane sections of projective
surfaces.\ {\it Manuscripta Math.\ \bf 68}, (1990) 17-34.
\vskip .2cm
\item{[Fl]} Fl{\o}ystad,G.:\ Determining obstructions for space curves,
with applications to non-reduced components of the Hilbert scheme.\
{\it J.\ Reine Angew.\ Math.\ \bf 439}, (1993), 11-44.
\vskip .2cm
\item{[Fu]} Fujita,T.:\ On the structure of polarized manifolds with
total deficiency one II.\ {\it J.\ Math.\ Soc.\ Japan \bf 33-3}, (1981)
415-434.
\vskip .2cm
\item{[G]} Green,M.:\ Koszul cohomology and
the geometry of projective varieties.\ {\it J.\ Diff.\ Geom.\ \bf 19},
(1984) 125-171.
\vskip .2cm
\item{[GP]} Gruson,L., Peskine,C.:\ Genre des courbes dans l'espace
projectif (II).\ {\it Ann.\ Sci.\ \`Ec.\ Norm.\ Sup.\ \bf 15}, (1982),
32-59.
\vskip .2cm
\item{[Ho]} Horowitz,T.:\ Varieties of low degree.\ {\it Brown
University Ph.\ D.\ Thesis}, (1982).
\vskip .2cm
\item{[Hi]} Hironaka,H.:\ On the arithmetic genera and the effective
genera of algebraic curves.\ {\it Mem.\ of the Coll.\ of Sci., Univ.\
of Kyoto}, ser.\ A, {\bf 30} (1957), 177-195.
\vskip .2cm
\item{[I]} Ionescu,P.:\ Vari\'et\'es projectives lisses de degr\'es 5
et 6.\ {\it C.\ R.\ Acad.\ Sci.\ Paris.\ \bf 293}, (1981) 685-687.
\vskip .2cm
\item{[K]} Kleppe,J.O.:\ Non-reduced components of the Hilbert scheme
of smooth space curves.\ In: {\it Proceedings Rocca di Papa 1982.\
Lecture Notes in Math.\ \bf 1266}.\ Springer, Berlin-New York:\ 1987,
180-207.
\vskip .2cm
\item{[L]} Lopez,A.F.: Surjectivity of Gaussian maps on curves in
$\I^r$ with general moduli. {\it Preprint}.
\vskip .2cm
\item{[M]} Mumford,D.:\ Further pathologies in algebraic geometry.\
{\it Amer.\ J.\ of Math.\ \bf 84}, (1962), 642-648.
\vskip .2cm
\item{[P]} Pinkham,H.: Deformations of algebraic varieties with $G_{m}$
action.\ {\it Harvard University Ph.\ D.\ Thesis}, (1974).
\vskip .2cm
\item{[S1]} Schreyer,F.O.:\ A standard basis approach to syzygies of
canonical curves.\  {\it J.\ Reine Angew.\ Math.\ \bf 421}, (1991)
83-123.   \vskip .2cm
\item{[S2]} Schreyer,F.O.:\ Syzygies of canonical curves and special
linear series.\  {\it Math.\ Ann.\ \bf 275}, (1986) 105-137.
\vskip .2cm
\item{[Sc]} Scorza,G.:\ Le variet\`a a curve sezioni ellittiche.\
{\it Ann.\ di Mat.\ }(3) {\bf 15}, (1908).
\vskip .2cm
\item{[Se1]} Serrano,F.:\ The adjunction mapping and hyperelliptic
divisors on a surface.\ {\it J.\ Reine Angew.\ Math.\ \bf 381}, (1987)
90-109.
\vskip .2cm
\item{[Se2]} Serrano,F.:\ Surfaces having a hyperplane section with a
special pencil.\ {\it Preprint}.
\vskip .2cm
\item{[St]} Stevens,J.:\ Deformations of cones over hyperelliptic
curves.\ {\it Preprint}.
\vskip .2cm
\item{[T]} Tendian,S.:\ Deformations of cones over curves of high
degree.\ {\it Ph. D. Thesis}, Univ.\ of North Carolina, 1990.
\vskip .2cm
\item{[V]} Voisin,C.:\ Courbes tetragonales et cohomologie de
Koszul.\  {\it J.\ Reine Angew.\ Math.\ \bf 387}, (1988) 111-121.
\vskip .2cm
\item{[W1]} Wahl,J.:\ The Jacobian algebra of a graded Gorenstein
singularity.\ {\it Duke Math.\ J.\ \bf 55}, (1987) 843-871.
\vskip .2cm
\item{[W2]} Wahl,J.:\ Introduction to Gaussian maps on an algebraic
curve.\ In:\ {\it Complex Projective Geometry, Trieste-Bergen 1989.
London Math.\ Soc.\ Lecture Notes Series \ \bf 179}.\ Cambridge Univ.\
Press:\  1992, 304-323.
\vskip .2cm
\item{[Z]} Zak,F.L.:\ Some properties of dual varieties and their
application in projective geometry.\ In:\ {\it Algebraic Geometry,
Proceedings Chicago 1989. \ Lecture Notes in Math.\ \bf 1479}.\
Springer, Berlin-New York:\ 1991, 273-280.

\vskip .5cm

ADDRESSES OF THE AUTHORS:
\vskip .2cm
Ciro Ciliberto, Dipartimento di Matematica, Universit\`a di Roma II,
Tor Vergata,

Viale della Ricerca Scientifica, 00133 Roma, Italy

e-mail: ciliberto@mat.utovrm.it
\vskip .2cm
Angelo Felice Lopez, Dipartimento di Matematica, Terza Universit\`a di
Roma,

Via Corrado Segre 2, 00146 Roma, Italy

e-mail: lopez@matrm3.mat.uniroma3.it
\vskip .2cm
Rick Miranda, Department of Mathematics, Colorado State University,
Ft. Collins,

CO 80523, USA

e-mail: miranda@math.colostate.edu

\end